\documentclass[11pt]{article}%
\usepackage{amsmath}
\usepackage{amsfonts}
\usepackage{amssymb}
\usepackage{graphicx}
\usepackage{a4wide}
\providecommand{\U}[1]{\protect\rule{.1in}{.1in}}
\newtheorem{theorem}{Theorem}

\newtheorem{corollary}[theorem]{Corollary}

\newtheorem{definition}[theorem]{Definition}
\newtheorem{example}[theorem]{Example}

\newtheorem{proposition}[theorem]{Proposition}

\begin{document}

\title{On Berinde's method for comparing iterative processes}
\author{Constantin Z\u{a}linescu\thanks{Octav Mayer Institute of Mathematics, Ia\c{s}i
Branch of Romanian Academy, Ia\c{s}i, Romania, email:
\texttt{zalinesc@uaic.ro}.}}
\date{}
\maketitle

\begin{abstract}
In the literature there are several methods for comparing two convergent
iterative processes for the same problem. In this note we have in view mostly
the one introduced by Berinde in [Picard iteration converges faster than Mann
iteration for a class of quasi-contractive operators, Fixed Point Theory and
Applications 2, 97--105 (2004)] because it seems to be very successful. In
fact, if IP1 and IP2 are two iterative processes converging to the same
element, then IP1 is faster than IP2 in the sense of Berinde. The aim of this
note is to prove this almost obvious assertion and to discuss briefly several
papers that cite the mentioned Berinde's paper and use his method for
comparing iterative processes.



\end{abstract}

\section{Introduction}

In the literature there are several methods for comparing two convergent
iterative processes for the same problem. In this note we have in view mostly
the one introduced by Berinde in \cite[Definition 2.7]{Ber04} because it seems
to be very successful. This was pointed out by Berinde himself in
\cite{Ber16}: \textquotedblleft This concept turned out to be a very useful
and versatile tool in studying the fixed point iterative schemes and hence
various authors have used it". However, it was pointed out by Popescu, using
\cite[Example 3.4]{Pop07}, that Berinde's method is not consistent. The
inconsistency of Berinde's method is indirectly mentioned also by Qing~\&
Rhoades in \cite[page 2]{QinRho08} by providing a very simple counterexample
in $\mathbb{R}$ to \cite[Theorem 2.1]{BabVar06}\footnote{Note that Berinde's
paper \cite{Ber04} is not cited in \cite{QinRho08}; see also \cite[Remark
3.2]{Ber16}.}. Moreover, referring to Berinde's method, Phuengrattana~\&
Suantai say in \cite[page 218]{PhuSua13}: \textquotedblleft It seem not to be
clear if we use above definition for comparing the rate of convergence". In
fact, if IP1 and IP2 are two (arbitrary) iterative processes converging to the
same element, then IP1 is faster than IP2 (and vice-versa) in the sense of
Berinde (\cite[Definition 2.7]{Ber04}).

The aim of this note is to prove this almost obvious assertion and to discuss
briefly several papers that cite \cite{Ber04} and refer to Berinde's method
for comparing iterative processes.

\section{Definitions and the main assertion\label{sec1}}

First, we quote from \cite[pages 99, 100]{Ber04} the text containing the
definitions which we have in view; these are reproduced in many papers from
our bibliography.

\medskip``Definition 2.5. Let $\{a_{n}\}_{n=0}^{\infty}$, $\{b_{n}%
\}_{n=0}^{\infty}$ be two sequences of real numbers that converge to $a$ and
$b$, respectively, and assume that there exists $l=\lim_{n\rightarrow\infty
}\big\vert
\frac{a_{n}-a}{b_{n}-b}\big\vert $.

(a) If $l=0$, then it can be said that $\{a_{n}\}_{n=0}^{\infty}$ converges
\emph{faster} to $a$ than $\{b_{n}\}_{n=0}^{\infty}$ to $b$.

(b) If $0<l<\infty$, then it can be said that $\{a_{n}\}_{n=0}^{\infty}$, and
$\{b_{n}\}_{n=0}^{\infty}$ \emph{have the same rate of convergence}."

\medskip

``Suppose that for two fixed point iteration procedures $\{u_{n}%
\}_{n=0}^{\infty}$ and $\{v_{n}\}_{n=0}^{\infty}$, both converging to the same
fixed point $p$, the error estimates

$\left\Vert u_{n}-p\right\Vert \leq a_{n}$, $n=0,1,2,..$. (2.7)

$\left\Vert v_{n}-p\right\Vert \leq b_{n}$, $n=0,1,2,..$. (2.8)

\noindent are available, where $\{a_{n}\}_{n=0}^{\infty}$ and $\{b_{n}
\}_{n=0}^{\infty}$ are two sequences of positive numbers (converging to zero).

Then, in view of Definition 2.5, we will adopt the following concept.

\medskip

Definition 2.7. Let $\{u_{n}\}_{n=0}^{\infty}$ and $\{v_{n}\}_{n=0}^{\infty}$
be two fixed point iteration procedures that converge to the same fixed point
$p$ and satisfy (2.7) and (2.8), respectively. If $\{a_{n}\}_{n=0}^{\infty}$
converges faster than $\{b_{n}\}_{n=0}^{\infty}$, then it can be said that
$\{u_{n}\}_{n=0}^{\infty}$ \emph{converges faster} than $\{v_{n}%
\}_{n=0}^{\infty}$ to $p$."

\medskip

Practically, the text above is reproduced in \cite[pages 30, 31]{Ber16},
getting so Definitions 1.1 and 1.2. The only differences are: ``(2.7)" and
``(2.8) are available, where" are replaced by ``(1.7)" and ``(1.8) are
available \emph{(and these estimates are the best ones available)}, where", respectively.

Immediately after \cite[Definition 1.2]{Ber16} it is said:

\medskip``This concept turned out to be a very useful and versatile tool in
studying the fixed point iterative schemes and hence various authors have used
it, see [1]-[5], [18], [22], [23], [28], [32]-[34], [37]-[41], [40],
[43]-[46], [55]-[57], [66], [68]-[72], [74], [78]-[81], to cite just an
incomplete list."\footnote{Throughout this paper the references mentioned in
the quoted texts are those in the works from where the texts are taken.}

\medskip

Note that Definition 9.1 from \cite{Ber07} is equivalent to Definition 2.5
from \cite{Ber04}; replacing $u_{n}$, $v_{n}$, $p$, $\left\Vert u_{n}
-p\right\Vert $ and $\left\Vert v_{n}-p\right\Vert $ with $x_{n}$, $y_{n}$,
$x^{\ast}$, $d(x_{n},x^{\ast})$ and $d(y_{n},x^{\ast})$ in (2.7), (2.8) and
Definition 2.7 from \cite{Ber04}, one obtains relations (5), (6) from
\cite[page 201]{Ber07} and an equivalent formulation of \cite[Definition
9.2]{Ber07}, respectively. Note that these definitions from Berinde's book
\cite{Ber07} are presented in the lecture \cite{Ber07b}.

\medskip

The proof of \cite[Theorem 3.1]{Ber04} is an example of the use of
\cite[Definition 2.7]{Ber04}; the text below is quoted from \cite{Ber04}, and
is followed by a short discussion:

\smallskip

\textquotedblleft The main result of this paper...

Theorem 3.1. Let $E$ be a uniformly convex Banach space, $K$ a closed convex
subset of $E$, and $T:K\rightarrow K$ a Zamfirescu operator; that is, an
operator that satisfies (z1), (z2), and (z3). Let $\{x_{n}\}_{n=0}^{\infty}$
be the Picard iteration associated with $T$, starting from $x_{0}\in K$, given
by (2.3), and $\{y_{n}\}$ the Mann iteration given by (2.1), where
$\{\alpha_{n}\}_{n=0}^{\infty}$ is a sequence satisfying (i) $\alpha_{1}=1$;
(ii) $0\leq\alpha_{n}<1$ for $n\geq1$; (iii) $\sum_{n=0}^{\infty}\alpha
_{n}(1-\alpha_{n})=\infty$. Then,

(1) $T$ has a unique fixed point in $E$, that is, $F_{T}=\{p\}$; (2) the
Picard iteration $\{x_{n}\}$ converges to $p$ for any $x_{0}\in K$; (3) the
Mann iteration $\{y_{n}\}$ converges to $p$ for any $y_{0}\in K$ and
$\{\alpha_{n}\}$ satisfying (i), (ii), and (iii); (4) Picard iteration is
faster than any Mann iteration.

Proof. Conclusions (1), (2), and (3) follow by Theorems 2.3 and 2.4. (4) ...
in view of the assumptions $0\leq a<1$; $0\leq b<1/2$; $0\leq c<1/2$ it
follows that $0\leq\delta<1$ ... which inductively yields

$\left\Vert x_{n+1}-p\right\Vert \leq\delta^{n}\left\Vert x_{1}-p\right\Vert $
(3.8)... $\left\Vert y_{n+1}-p\right\Vert \leq{\textstyle\prod\nolimits_{k=1}%
^{n}} \left[  1-\alpha_{k}+3\delta\alpha_{k}\right]  \cdot\left\Vert
y_{1}-p\right\Vert $ (3.12) ...

\noindent if $\delta\in\lbrack0,1/3)$ then $0<1-\alpha_{k}+3\delta\alpha
_{k}<1$, (3.13) while for $\delta\in\lbrack1/3,1)$ we have $1-\alpha
_{k}+3\delta\alpha_{k}\geq1$. (3.14) Thus, for $\delta\in\lbrack1/3,1)$ we
have $0\leq\lim_{n\rightarrow\infty}\frac{\delta^{n}}{{\textstyle\prod
\nolimits_{k=1}^{n}} (1-\alpha_{k}+3\delta\alpha_{k})}\leq\lim_{n\rightarrow
\infty}\delta^{n}=0$ (3.15) ... If $\delta\in\lbrack0,1/3)$ ... which yields
$\frac{\delta}{1-\alpha_{k}+3\delta\alpha_{k}}<(1-\delta)$ (3.17) ... Hence
$\lim_{n\rightarrow\infty}\frac{\delta^{n}}{{\textstyle\prod\nolimits_{k=1}%
^{n}} (1-\alpha_{k}+3\delta\alpha_{k})}=0$. \ (3.19)"

\medskip Hence, in both cases, even not said explicitly, $a_{n}:=\delta
^{n}\left\Vert x_{1}-p\right\Vert $ and $b_{n}:={\textstyle\prod
\nolimits_{k=1}^{n}} \left[  1-\alpha_{k}+3\delta\alpha_{k}\right]
\cdot\left\Vert y_{1}-p\right\Vert $; clearly $a_{n}\rightarrow0$ in both
cases, while $b_{n}\rightarrow0$ if $\delta\in\lbrack0,1/3)$, and so (4) holds
by \cite[Definition 2.7]{Ber04} in this case.

Assume that $\delta\in\lbrack1/3,1)$. If $\delta=1/3$ or $y_{1}=p$ then
$y_{n}=p$ for $n\geq1$ and so $\{y_{n}\}$ converge faster than $\{x_{n}\}$.
So, take $\delta\in(1/3,1)$ and $y_{1}\neq p$; then $b_{n}\geq\left\Vert
y_{1}-p\right\Vert =:\beta>0$ for every $n\geq1$ and so $\{x_{n}\}$ converges
faster than $\{y_{n}\}$ by \cite[Definition 2.7]{Ber04} if one does not take
into account \textquotedblleft(\emph{converging to zero})" from the preamble
of \cite[Definition 2.7]{Ber04}. In fact, by (3.14), the sequence $\{b_{n}\}$
is increasing and so there exists $b\in\lbrack\beta,\infty]$ such that
$b_{n}\rightarrow b$; moreover, $b\in\mathbb{R}$ if and only if the series
$\sum_{n=1}^{\infty}\alpha_{n}$ is convergent.

\medskip Possibly, the proof of \cite[Theorem 3.1]{Ber04} represents a
motivation for not asking the convergence to $0$ of $(a_{n})$ and $(b_{n})$ in
the preamble of \cite[Definition 2.7]{Ber04} (and \cite[Definition 9.2]%
{Ber07}, \cite[Definition 1.2]{Ber16}), and also for the absence of
\textquotedblleft(\emph{converging to zero})" in \cite[Def.\ 2]{AbAsDe22},
\cite[page 8]{AgSaIs25}, \cite[Def.\ 1.3]{AlGuAk18}, \cite[Def.\ 2]{AlaRoh25},
\cite[Def.\ 2.2]{AlRoTo25}, \cite[Def.\ 1.3]{AliAli20b}, \cite[Def.\ 1.10]%
{AliAli20}, \cite[Def.\ 5]{AlAlDaTuZaMa23}, \cite[Def.\ 1.3]{AlAlRo21},
\cite[Def.\ 1.4]{AlAlUd21}, \cite[Def.\ 4]{AlHuKaCh22}, \cite[Def.\ 1.3]%
{AlJuAl22}, \cite[Def.\ 3]{AnsMas23}, \cite[Def.\ 1.1]{BhuTiw19},
\cite[Def.\ 2.10]{DunHie20},\footnote{Notice that in this paper one recalls
(as Definition 2.12) and uses only \cite[Definition 3.5]{Pop07}.} \cite[page
3]{FaGhPoRe15}, \cite[page 1518]{FatRez18}, \cite[Def.\ 2.4]{HaOuChMo24},
\cite[Def.\ 5]{HuHuAl21}, \cite[Def.\ 2.5]{HuUlAr18}, \cite[Def.\ 3]%
{JuAlKu22}, \cite[Def.\ 3]{KauCha22b}, \cite[Def.\ 2.8]{PiDaRaGh19},
\cite[Def.\ 2.8]{RaReDaGh22}. Notice that $(b_{n})$ or/and $(a_{n})$ are
constant sequences in \cite{BabVar06, BabVar07}, \cite{Kum14},
\cite{FaGhPoRe15}, \cite{VeJaSh16}, \cite{FatRez18}, \cite{AkhKha19},
\cite[Proposition 1]{OkAbDe20},\footnote{This is a special situation having in
view its conclusion, $X$ being \textquotedblleft a normed linear
space\textquotedblleft\ and \textquotedblleft$x_{0}=u_{0}\in C$". One lets
\textquotedblleft$\left\vert a_{n}-x^{\ast}\right\vert =\left\Vert
x_{0}-x^{\ast}\right\Vert $. (45) ... $\left\vert b_{n}-x^{\ast}\right\vert
=\left\Vert u_{0}-x^{\ast}\right\Vert $. (51) Hence, using (45), (51), and the
condition that $x_{0}=u_{0}\in C$, we obtain
\[
\lim_{n\rightarrow\infty}\frac{\left\vert a_{n}-x^{\ast}\right\vert
}{\left\vert b_{n}-x^{\ast}\right\vert }=\lim_{n\rightarrow\infty}%
\frac{\left\Vert x_{0}-x^{\ast}\right\Vert }{\left\Vert u_{0}-x^{\ast
}\right\Vert }=\frac{\left\Vert x_{0}-x^{\ast}\right\Vert }{\left\Vert
u_{0}-x^{\ast}\right\Vert }=\frac{\left\Vert x_{0}-x^{\ast}\right\Vert
}{\left\Vert x_{0}-x^{\ast}\right\Vert }=1.
\]
Because $0<l=1<\infty$, it follows that the sequences $\{x_{n}\}$ and
$\{u_{n}\}$ have the same rate of convergence", even if in \cite[Definition
6]{OkAbDe20} one mentions \textquotedblleft where $\{a_{n}\}_{n=0}^{\infty}$
and $\{b_{n}\}_{n=0}^{\infty}$ are two sequences of positive numbers
converging to zero".} \cite{JuAlKu22}, \cite{KauCha22}, \cite{Ahm24},
\cite{AgSaIs25}; even more, like in the proof of \cite[Theorem 3.1]{Ber04},
one has that $b_{n}\rightarrow\infty$ in several papers.

\medskip

In the next result we use the version for metric spaces of \cite[Definition
2.7]{Ber04} (see \cite[Definition 9.2]{Ber07}).

\begin{proposition}
\label{p-fast}Let $(X,d)$ be a metric space and $(x_{n})_{n\geq1}$,
$(y_{n})_{n\geq1}$ be two sequences from $X$ converging to $x^{\ast}\in X$.
Then $(x_{n})$ converges faster than $(y_{n})$ to $x^{\ast}$.
\end{proposition}

Proof. For each $n\geq1$ let us consider
\[
0<a_{n}:=d(x_{n},x^{\ast})+d(y_{n},x^{\ast})+\frac{1}{n},\quad0<b_{n}
:=\left\{
\begin{array}
[c]{ll}%
\sqrt{a_{n}} & \text{if }a_{n}\leq1,\\
d(y_{n},x^{\ast}) & \text{otherwise.}%
\end{array}
\right.
\]
It follows that $a_{n}\rightarrow0$, $b_{n}\rightarrow0$,
\[
d(x_{n},x^{\ast})\leq a_{n},\quad d(y_{n},x^{\ast})\leq b_{n},\quad\forall
n\geq1,
\]
and $a_{n}/b_{n}=\sqrt{a_{n}}$ for sufficiently large $n;$ it follows that
$\lim_{n\rightarrow\infty}a_{n}/b_{n}=\lim_{n\rightarrow\infty}\sqrt{a_{n}}
=0$. Therefore, by \cite[Definition 9.2]{Ber07}, $(x_{n})$ converges faster to
$x^{\ast}$ than $(y_{n})$ does. \hfill$\square$

\begin{corollary}
\label{c-fast} Let $(X,d)$ be a metric space and $(x_{n})_{n\geq1}\subset X$
converging to $x^{\ast}\in X$ and set $y_{n}:=x^{\ast}$ for $n\geq1$. Then
$(x_{n})$ converges faster than $(y_{n})$ to $x^{\ast}$.
\end{corollary}

From our point of view, the preceding \textquotedblleft
result\textquotedblright\ shows that Berinde's notion of rapidity of
convergence for fixed point iterative schemes, recalled above, is not useful,
even if Berinde in \cite[page 35]{Ber16} claims that \textquotedblleft Of all
concepts of rapidity of convergence presented above for numerical sequences,
the one introduced by us in Definition 1.2 [14] appears to be the most
suitable in the study of fixed point iterative methods". Berinde (see
\cite[page 36]{Ber16}) mentions that he \textquotedblleft tacitly admitted in
Definition 1.2 that \emph{the estimates (1.7) and (1.8)} taken into
consideration \emph{are the best possible}". Clearly, \textquotedblleft the
estimates are the best ones available" and \textquotedblleft the estimates~...
are the best possible" are very different in meaning.\footnote{\label{note-3}%
Among the 35 papers from our bibliography published in the period 2017--2020,
our reference \cite{Ber16} is mentioned only in \cite{ErtGur19},
\cite{GuEkKhKa19}, \cite{GuErAb20} and \cite{GuKhErKa20}. However,
\cite[Definition 2.7]{Ber04} is used in \cite{ErtGur19}, \cite{GuEkKhKa19} and
\cite{GuKhErKa20} without any mention that the obtained estimates are the best
possible.}

\smallskip Of course, \emph{the best possible estimates in relations (1.7) and
(1.8)} from \cite{Ber16} (that is in relations (2.7) and (2.8) from
\cite{Ber04} recalled above) \emph{are}
\begin{equation}
a_{n}:=\left\Vert u_{n}-p\right\Vert ,\quad b_{n}:=\left\Vert v_{n}%
-p\right\Vert \quad(n\geq0). \label{r-1}%
\end{equation}

Assuming that $d(x_{n},x^{\ast})\rightarrow0$, getting (better) upper
estimates for $d(x_{n},x^{\ast})$ depends on the proof, including the author's
ability to majorize certain expressions. Surely, \emph{the best available
estimates are} exactly \emph{those obtained by the authors in their proofs}
(when not using estimates obtained by other authors as in \cite{AbbNaz14}
(\cite{Sah11}), \cite{Mog14}~\& \cite{Mog16} (\cite{Kha13}), \cite{ThThPo14}%
~\& \cite{ThThPo16b} (\cite{AbbNaz14}), \cite{Gur15}
(\cite{HuaNoo07}), \cite{AsKhAl16} (\cite{ThThPo16b}), \cite{Kha16}
(\cite{ThThPo16b}), \cite{SaPaTi16} (\cite{AbbNaz14}),
\cite{VeJaSh16} (\cite{SaAnYa15}), \cite{AlGuAk20}
(\cite{KadYil14}~\& \cite{ThThPo16b}), \cite{AliAli20}
(\cite{Kha13}~\& \cite{SinPit16}), \cite{AliAli20b} (\cite{Sah11}~\&
\cite{UllArs16}), \cite{AlAlRo21} (\cite{Kha13}~\&  \cite{Mog16}),
\cite{KaRaSuBiGhInAl21} (\cite{ThThPo16b}), \cite{AbAsDe22} (\cite{AbbNaz14}%
~\& \cite{ThThPo16b}), \cite{AlHuKaCh22} (\cite{HuUlAr18}~\& \cite{UllArs18b}%
); \cite{AlJuAl22} (\cite{GurKar14}); \cite{JiShAhShBo22}
(\cite{Sri22}), \cite{KauCha22b} (\cite{AbbNaz14}, \cite{ThThPo16},
\cite{SinPit16}, \cite{UllArs18b}~\& \cite{AliAli20}), \cite{Kha22}
(\cite{Sah11}, \cite{GurKar14}~\& \cite{UllArs16}),
\cite{AlAlDaTuZaMa23} (\cite{Gur16}), \cite{AnsMas23}
(\cite{AliAli20}, \cite{AliAli20b}~\& \cite{HuHuAl21}),
\cite{HamKat23} (\cite{UllArs18}), \cite{Beg24} (\cite{ThThPo16b}~\&
\cite{RaReDaGh22}), \cite{GauVin24} (\cite{Kha13},
\cite{SinPit16}~\& \cite{AliAli20}), \cite{RaReDaGh24}
(\cite{ThThPo16b}).

\smallskip The use of Berinde's method for comparing the speeds of convergence
is very subjective. It is analogue to deciding that $a/b\leq c/d$ knowing only
that $0<a\leq c$ and $0<b\leq d$\thinspace! Effectively, there are several
works in which one writes directly that $\lim_{n\rightarrow\infty}%
\frac{\left\Vert u_{n}-p\right\Vert }{\left\Vert v_{n}-p\right\Vert }\leq
\lim_{n\rightarrow\infty}\frac{a_{n}}{b_{n}}$ if one got $\left\Vert
u_{n}-p\right\Vert \leq a_{n}$ and $\left\Vert v_{n}-p\right\Vert \leq b_{n}$
for $n\geq1$ (see \cite{PanMis23}), or even $\lim_{n\rightarrow\infty}%
\frac{\left\Vert u_{n}-p\right\Vert }{\left\Vert v_{n}-p\right\Vert }%
=\lim_{n\rightarrow\infty}\frac{a_{n}}{b_{n}}$ (see
\cite{KaBoDoAt15}, \cite{KaAtDoBo16}, \cite{UllArs16},
\cite{UdOfIg21b}, \cite{HaReDe22}, \cite{SaBaGu22}, \cite{SahBan23},
\cite{SalAbe23}, \cite{SahBan24}).

\smallskip Taking $a_{n}$ and $b_{n}$ defined by (\ref{r-1}) in
\cite[Definition 2.7]{Ber04} one obtains Definition 3.5 of Popescu from
\cite{Pop07}\footnote{\label{note-x=r}Of course, when $(X,\left\Vert
\cdot\right\Vert )$ is $(\mathbb{R},\left\vert \cdot\right\vert )$, Definition
3.5 of Popescu \cite{Pop07} reduces to Definition 2.5(a) of Berinde
\cite{Ber04} when $a=b$.}. Popescu's definition is used explicitly by
Rhodes~\& Xue (see \cite[page 3]{RhoXue10}), but they wrongly atribute it to
\cite{Ber04}; this attribution is wrong because \cite[Definition 3.5]{Pop07}
reduces to \cite[Definition 2.5]{Ber04} only in the case in which the involved
normed vector space is $\mathbb{R}$. Note that Rhoades knew about Popescu's
definition because \cite{Pop07} is cited in \cite[page 2]{QinRho08}.

\smallskip Notice that Popescu's definition is extended to metric spaces by
Berinde, Khan~\& P\u{a}curar in \cite[page 8]{BeKhPa15}, as well as by
Fukhar-ud-din~\& Berinde in \cite[page 228]{FukBer16}; also observe that
Popescu's paper \cite{Pop07} is not cited in \cite{BeKhPa15} and
\cite{FukBer16}.

\smallskip Even if in \cite{Ber04} it is not defined when two iteration
schemes have the same rate of convergence, Dogan~\& Karakaya obtain that
\textquotedblleft the iteration schemes $\{k_{n}\}_{n=0}^{\infty}$ and
$\{l_{n}\}_{n=0}^{\infty}$ have the same rate of convergence to $p$ of $\wp$"
in \cite[Theorem 2.4]{DogKar18}; the conclusion of \cite[Theorem
2.4]{DogKar18} is based on the fact that its authors found the same upper
estimates for $\left\Vert k_{n+1}-p\right\Vert $ and $\left\Vert
l_{n+1}-p\right\Vert $ when $l_{0}=k_{0}$.

Accepting such an argument, and taking $a_{n}:=b_{n}:=d(x_{n},x^{\ast
})+d(y_{n},x^{\ast})+\frac{1}{n}$ in the proof of Proposition \ref{p-fast},
one should obtain that any pair of sequences $(x_{n})_{n\geq1}$,
$(y_{n})_{n\geq1}\subset(X,d)$ with the same limit $x^{\ast}\in X$ have the
same rate of convergence.

\smallskip Recall that Rhoades in \cite[pages 742, 743]{Rho76} says that
having \textquotedblleft$\{x_{n}\}$, $\{z_{n}\}$ two iteration schemes which
converge to the same fixed point $p$, we shall say that $\{x_{n}\}$ is better
than $\{z_{n}\}$ if $\left\vert x_{n}-p\right\vert \leq\left\vert
z_{n}-p\right\vert $ for all $n$"; having in view the previous definition and
\cite[Example 2.8]{Ber04}, Berinde claims that \textquotedblleft The previous
example shows that Definition 2.7 introduces a sharper concept of rate of
convergence than the one considered by Rhoades [11]\textquotedblright. In this
context we propose the following definition.

\begin{definition}
\label{d-better}Let $(X,d)$ be a metric space, and let $(x_{n})_{n\geq1}$,
$(y_{n})_{n\geq1}\subset(X,d)$ and $x,y\in X$ be such that $x_{n}\rightarrow
x$, $y_{n}\rightarrow y$. One says that $(x_{n})$ converges better to $x$ than
$(y_{n})$ to $y$ if there exists some $\alpha>0$ such that $d(x_{n}
,x)\leq\alpha d(y_{n},y)$ for sufficiently large $n;$ one says that $(x_{n})$
and $(y_{n})$ have the same rate of convergence if $(x_{n})$ converges better
to $x$ than $(y_{n})$ to $y$, and $(y_{n})$ converges better to $y$ than
$(x_{n})$ to $x$.
\end{definition}

Using the conventions $\frac{0}{0}:=1$ and $\frac{\alpha}{0}:=\infty$ for
$\alpha>0$, [$(x_{n})$ converges better to $x$ than $(y_{n})$ to $y$] if and
only if $\limsup_{n\rightarrow\infty}\frac{d(x_{n},x)}{d(y_{n},y)}<\infty$;
consequently, [$(x_{n})$ and $(y_{n})$ have the same rate of convergence] (in
the sense of Definition \ref{d-better}) if and only if $0<\liminf
_{n\rightarrow\infty}\frac{d(x_{n},x)}{d(y_{n},y)}\leq\limsup_{n\rightarrow
\infty}\frac{d(x_{n},x)}{d(y_{n},y)}<\infty$.\footnote{It seems that this way
of comparing the rate of convergence for sequences of real numbers was
introduced by Knopp in \cite{Kno47} (cf.\ \cite[Definition 1.2]{GuKhFu17}).}

\begin{example}
\label{ex1}Consider the sequences $(x_{n})_{n\geq1}$, $(y_{n})_{n\geq1}%
\subset\mathbb{R}$ defined by
\[
x_{n}:=\left\{
\begin{array}
[c]{ll}%
n^{-1} & \text{if }n\text{ is odd,}\\
(2n)^{-1} & \text{if }n\text{ is even,}%
\end{array}
\right.  \quad y_{n}:=\left\{
\begin{array}
[c]{ll}%
(2n)^{-1} & \text{if }n\text{ is odd,}\\
n^{-1} & \text{if }n\text{ is even.}%
\end{array}
\right.
\]
Clearly $\lim_{n\rightarrow\infty}x_{n}=\lim_{n\rightarrow\infty}y_{n}=0$, and
it is very natural to consider that they have the same rate of convergence;
this is confirmed using Definition \ref{d-better}. It is obvious that neither
$(x_{n})$ is better (faster) than $(y_{n})$, nor $(y_{n})$ is better (faster)
than $(x_{n})$ in the senses of Rhoades (\cite{Rho76}), or Berinde
\cite{Ber04}, or Popescu \cite{Pop07}, or Berinde, Khan~\&
P\u{a}curar\ (\cite{BeKhPa15}), or Fukhar-ud-din~\& Berinde (\cite{FukBer16}).
\end{example}

It is interesting that Definitions 2.5\thinspace(a) and 2.7 from \cite{Ber04}
(recalled above) are transformed into results (that is, true logical
propositions) in some articles. The following text is quoted from
\cite[p.\ 302]{MebMew19}, where \textquotedblleft\lbrack6]" is our reference
\cite{Ber04}; replacing [6] by [7] one gets a text from \cite[p.\ 2300]%
{AkNaAfMe21}; see also \cite[Lemma 2.8]{SaBaGu22}, \cite[Lemma 2.5]%
{MePiNaOnAd23}, \cite[Lemma 2.10]{SahBan23}, \cite[Lemma 1]{SahBan24},
\cite[Lemma 2]{ShRaKaBeRa25}.

\medskip

\textquotedblleft Lemma 2.1. [6] Let $\left\{  a_{n}\right\}  $ and $\left\{
b_{n}\right\}  $ be two sequences of real numbers converging to $a$ and $b$
respectively. If $\lim_{n\rightarrow\infty}\frac{\left\vert a_{n}-a\right\vert
}{\left\vert b_{n}-b\right\vert }=0$, then $\{a_{n}\}$ converges faster than
$\left\{  b_{n}\right\}  $.

\smallskip Lemma 2.2. [6] Suppose that for two fixed point iteration processes
$\left\{  u_{n}\right\}  $ and $\left\{  v_{n}\right\}  $ both converging to
the same fixed point $x^{\ast}$, the error estimates $\left\Vert u_{n}%
-x^{\ast}\right\Vert \leq a_{n}$ \ $\ n\geq1$, $\left\Vert v_{n}-u^{\ast
}\right\Vert \leq b_{n}$ \ $n\geq1$, are available where $\left\{
a_{n}\right\}  $ and $\left\{  b_{n}\right\}  $ are two sequences of positive
numbers converging to zero. If $\left\{  a_{n}\right\}  $ converges faster
than $\left\{  b_{n}\right\}  $, then $\left\{  u_{n}\right\}  $ converges
faster than $\left\{  v_{n}\right\}  $ to $x^{\ast}$."

\smallskip Let us comment the following text from \cite{AliAli20b}%
:\footnote{Replacing [14], Definition 1.3 and [5] from \cite[Remark
1.4]{AliAli20b} by [8], Definition 1.4 and [9] (respectively), one obtains
\cite[Remark 1.5]{AlAlUd21}.}

\smallskip\textquotedblleft\textbf{Remark 1.4} \emph{In 2007, Popescu [14]
claimed that Definition 1.3 is not consistent and gave a new
definition (see Definition 3.5 [14]) to compare the rate of
convergence of iterative methods which is almost same as Definition
1.3 where he just replaced \textquotedblleft$\leq$" sign with
\textquotedblleft$=$" sign. On careful reading of Popescu's paper,
we found that he also used Definition 1.3 in his result (Theorem 3.7
[14], see error bounds in (3.15), (3.19) and (3.22)) and even bounds
obtained in (3.15) and (3.22) are not equal as claimed by him. Quite
recently, Berinde [5] wrote a review paper and clarified that
Popescu's claim is not correct in general.}"\footnote{Notice that
the conclusions of \cite[Remark 1.4]{AliAli20b} are taken up by Saif
et~al.\ in \cite[page 440]{SaAlAlUd24}: \textquotedblleft Berinde
(2016) and Ali and Ali (2020a) found that Popesecu's claim is not
true in general," where Berinde (2016) and Ali and Ali (2020a) are
our references \cite{Ber16} and \cite{AliAli20b}, respectively.}

\medskip

1) \textquotedblleft[14]" and \textquotedblleft[5]" are our references
\cite{Pop07} and \cite{Ber16}, respectively, while \textquotedblleft
Definition 1.3" is equivalent to \cite[Definition 2.7]{Ber04} (less
``converging to zero").

\smallskip2) One says: \textquotedblleft... he also used Definition 1.3 in his
result (Theorem 3.7 [14]), see error bounds in (3.15), (3.19) and (3.22) and
even bounds obtained in (3.15) and (3.22) are not equal as claimed by him."
Probably \textquotedblleft not true" instead of \textquotedblleft not equal"
in the preceding text!

\smallskip3) \textquotedblleft Theorem 3.7 [14]" is: \textquotedblleft%
\emph{Let $E$ be an arbitrary Banach space, $K$ a closed convex subset of $E$,
and $T:K\rightarrow K$ a quasi-$\delta$-contraction. Let $\{y_{n}%
\}_{n=0}^{\infty}$ be defined by (3.1) and $y_{0}\in K$, $y_{0}\not \in F(T$)
with $\{\alpha_{n}\}\subset[0,1]$ satisfying}

\emph{(i)} $\sum\nolimits_{n=0}^{\infty}\alpha_{n}=\infty$.

\emph{Then }$\{y_{n}\}_{n=0}^{\infty}$\emph{ converges strongly to the fixed
point of }$T$\emph{ and, moreover, the Picard iteration }$\{x_{n}%
\}_{n=0}^{\infty}$\emph{ defined by (3.3) and }$x_{0}\in K$\emph{ converges
faster than the Mann iteration if}

\emph{(ii) }$\alpha_{n}<\frac{1}{1+\delta},\ \ n=0,1,2,$\emph{
...,\ \ \ \ \ (iii) }$\lim_{n\rightarrow\infty}\prod\nolimits_{k=0}^{n}\left[
\frac{\delta}{1-(1+\delta)\alpha_{k}}\right]  =0.$"

\smallskip4) Let us verify the inequality

\smallskip

\textquotedblleft$\left\Vert y_{n+1}-p\right\Vert \leq{\prod\nolimits_{k=0}%
^{n}}[1-(1-\delta)\alpha_{k}]\cdot\left\Vert y_{0}-p\right\Vert ,\ \ n=0,1,2,$
... (3.15)"

\smallskip First, $T$ is a \emph{quasi-}$\delta$\emph{-contraction }if
$\delta\in\lbrack0,1)$ and there exists $L>0$ such that (3.9) holds:

$\left\Vert Tx-Ty\right\Vert \leq\delta\left\Vert x-y\right\Vert +L\cdot
\min\{\left\Vert x-Tx\right\Vert ,\left\Vert y-Ty\right\Vert ,\left\Vert
x-Ty\right\Vert ,\left\Vert y-Tx\right\Vert \}$ $\forall x,y\in E$.

\smallskip

\textquotedblleft Proof. Using (3.1)" that is, \textquotedblleft%
$x_{n+1}=(1-\alpha_{n})x_{n}+\alpha_{n}Tx_{n},\ n=0,1,2,...$", and
$Tp=p$ \textquotedblleft we get" indeed \textquotedblleft(3.12) Take
$x:=p$ and $y:=y_{n}$ in (3.9) we obtain", indeed,
\textquotedblleft$\left\Vert Ty_{n}-p\right\Vert
\leq\delta\left\Vert y_{n}-p\right\Vert $, (3.13)", because
$\left\Vert Tx-x\right\Vert =0$ in (3.9), \textquotedblleft and
then", obviously, \textquotedblleft... (3.14)" by (3.12) and (3.13).
\textquotedblleft By induction, we get", clearly,
\textquotedblleft... (3.15)"; surely, it was preferable to write
$\left\Vert y_{0}-p\right\Vert $ before
${\prod\nolimits_{k=0}^{n}}$. So, the upper bound in (3.15) is
correct, contrary to what is said in \cite[Remark 1.4]{AliAli20b}.

\smallskip Using (i) one gets (3.17), that is, $y_{n}\rightarrow p$.

\smallskip5) In fact it is \textquotedblleft$\leq$" instead of
\textquotedblleft$<$" in (3.18) and so, instead of (3.19) one has

\smallskip$\left\Vert x_{n+1}-p\right\Vert \leq\delta^{n+1}\cdot\left\Vert
x_{0}-p\right\Vert ,\ \ n\geq0.$ \ \ (3.19')

\smallskip\noindent(whence $x_{n}\rightarrow p$). So, indeed, the inequality
in (3.19) could not be true! Set $a_{n}:=\delta^{n+1}\cdot\left\Vert
x_{0}-p\right\Vert $ for $n\geq1$.

\smallskip6) The aim is to get lower estimates for $\left\Vert y_{n}%
-p\right\Vert $ after (3.19), and so one had to mention that conditions (ii)
and (iii) hold, what we assume in the sequel; hence $[1-(1+\delta)\alpha
_{k}]>0$ by (ii). Then \textquotedblleft by (3.1) we have ... (3.20). Using
(3.13) we get ... (3.21) which implies that

\smallskip$\left\Vert y_{n+1}-p\right\Vert \geq{\prod\nolimits_{k=0}^{n}%
}[1-(1+\delta)\alpha_{k}]\left\Vert y_{0}-p\right\Vert ,\ \ n=0,1,2,$ .... (3.22)"

\smallskip Clearly $c_{n}:=\left\Vert y_{0}-p\right\Vert \cdot{\prod
\nolimits_{k=0}^{n}}[1-(1+\delta)\alpha_{k}]>0$. Using (3.19') and (3.22) one
obtains that
\[
(0\leq)\ \ \frac{\left\Vert x_{n+1}-p\right\Vert }{\left\Vert y_{n+1}%
-p\right\Vert }\leq\frac{a_{n}}{c_{n}}=\frac{\left\Vert x_{0}-p\right\Vert
}{\left\Vert y_{0}-p\right\Vert }\cdot{\prod\nolimits_{k=0}^{n}}\frac{\delta
}{[1-(1+\delta)\alpha_{k}]}\ \ \forall n\geq0,
\]
and so $\lim_{n\rightarrow\infty}\frac{\left\Vert x_{n}-p\right\Vert
}{\left\Vert y_{n}-p\right\Vert }=0$ by (iii).

\smallskip7) Consequently, Definition 1.3 was not used in the proof of
\cite[Theorem 3.7]{Pop07} and the only (typing?) errors are present in (3.18)
and (3.19) in which \textquotedblleft$<$" had to be replaced by
\textquotedblleft$\leq$".

\smallskip8) Related to \textquotedblleft Quite recently, Berinde [5] wrote a
review paper and clarified that Popescu's claim is not correct in general",
observe that the word \textquotedblleft Popescu" appears 3 times in
\cite{Ber16}: 2 times in the bibliography and once more in the text
\textquotedblleft Popescu [71], [72], compared Picard iteration and Mann
iteration in the class of so called quasi-$\varphi$-contractions, thus
extending significantly the results in [14], [16], [5] and [80]". Moreover,
\textquotedblleft\lbrack71]" appears 2 times (already mentioned) and
\textquotedblleft\lbrack72]\textquotedblright\ in other 4 places:

\textquotedblleft This concept turned out to be a very useful and versatile
tool... [66], [68]-[72], [74],...",

\textquotedblleft On the other hand, almost all the authors of the papers
...[66], [68]-[72], [74], ...",

\textquotedblleft As mentioned in Introduction, ... [66], [68]-[72], [78]-[81],...",

\textquotedblleft For other papers that used the concept ... [66], [68]-[72],
[74], ..."

\smallskip So, in \textquotedblleft Berinde [5]" we did not find any
assertion that \textquotedblleft clarified that Popescu's claim is
not correct in general"; even more, one says that Popescu extended
\textquotedblleft significantly the results in [14], [16], [5] and
[80]".

\section{Remarks on the use of Berinde and Popescu's notions in papers citing
\cite{Ber04}}

Practically, all the papers mentioned in the sequel were found on internet
when searching, with Google Scholar, the works citing Berinde's article
\cite{Ber04}.

\medskip First we give the list of articles, mentioning their authors and
results, in which Berinde's Definition 2.7 from \cite{Ber04} is used (even if
not said explicitly sometimes, or mentioning only \cite[Definition
2.5(a)]{Ber04}):

\medskip Berinde~\& Berinde~-- \cite[Theorem 3.3]{BerBer05}; Babu~\& Prasad~--
\cite[Theorems 3.1, 3.3]{BabVar06b} and \cite[Theorem 2.1]{BabVar06}
(+ \cite[Theorem 2.1]{BabVar07});\footnote{Notice the inequality
\textquotedblleft$a_{n}/b_{n}\leq a_{n}$" in the proof of Theorem
2.1!} Olaleru~-- \cite[Theorem 2]{Ola07b}, \cite[Theorem 1]{Ola07}
and \cite[Theorems 1, 2]{Ola09}\footnote{In fact, having in view
that $\delta:=k/(1-k)>1$ for $k\in(1/2,1)$, $\{a_{n}\}$ and
$\{b_{n}\}$ do not converge to $0$ in the results (at least in this
case) in these three papers; moreover, V Kumar (see \cite[page
1320]{Kum13}) shows that \cite[Theorem 2]{Ola09} is false by using
\cite[Definition 3.5]{Pop07} for a simple example in $\mathbb{R}$.};
Sahu~-- \cite[Theorem 3.6]{Sah11}; Akbulut~\& \"{O}zdemir~--
\cite[Theorem 2.3]{AkbOzd12}; Hussain et al.~-- \cite[Theorems 18,
19]{HuKuKu13}; Karahan~\& Ozdemir~-- \cite[Theorem 1]{KarOzd13};
Khan~-- \cite[Proposition 1]{Kha13}; Abbas~\& Nazir~-- \cite[Theorem
3]{AbbNaz14}; G\"{u}rsoy~\& Karakaya~-- \cite[Theorem 3]{GurKar14};
Kadioglu~\& Yildirim~-- \cite[Theorem 5]{KadYil14}; Karakaya et
al.~-- \cite[Theorem 3]{KaGuEr14} and \cite[Theorem 2.2]{KaGuEr16};
Kumar~-- \cite[Theorem 3.1]{Kum14}; Mogbademu~-- \cite[Theorem
3.1]{Mog14} and \cite[Theorem 2.1]{Mog16}; \"{O}zt\"{u}rk \c{C}eliker~--
\cite[Theorem 8]{Ozt14}; Thakur et al.~-- \cite[Theorem
2.3]{ThThPo14}\footnote{In \cite[page 3]{ThThPo14} one appreciates
that \textquotedblleft In recent years, Definition 2.2 has been used
as a standard tool to compare the fastness of two fixed point
iterations", Definition 2.2 being \cite[Definition 2.7]{Ber04}.} and
\cite[Theorem 3.1]{ThThPo16b}; Fathollahi et al.~--
\cite[Propositions 3.1, 3.2, Theorems 3.1, 4.1--4.4, Lemmas
3.1--3.4]{FaGhPoRe15}; G\"{u}rsoy~-- \cite[Theorem 3]{Gur15}; Jamil~\&
Abdullateef~-- \cite[Theorem 3.2]{JamAbd15}; Jamil~\& Abed~--
\cite[Theorems 3.1--3.4]{JamAbe15} and \cite[Theorems
3.1--3.4]{JamAbe15b}; Karakaya et al.~-- \cite[Theorem
5]{KaBoDoAt15}, \cite[Theorem 2.5]{KaAtDoBo16} and \cite[Theorem
2.4]{KaAtDoBo17}; Yadav~-- \cite[Example 2]{Yad15}; Abed~\& Abbas,
\cite[Theorem (3.8)]{AbeAbb16}; Asaduzzaman et al.~-- \cite[Theorem
3.3]{AsKhAl16}; Rani~\& Jyoti~--\cite[Theorem 13]{RanJyo16};
Khatun~-- \cite[Theorem 3.6.1]{Kha16}; Sahu et al.~-- \cite[Theorem
4.1]{SaPaTi16}; Sintunavarat~\& Pitea~-- \cite[Theorem
2.1]{SinPit16}; Ullah~\& Arshad~-- \cite[Theorem 4]{UllArs16}; Verma
et~al.~-- \cite{VeJaSh16};\footnote{See the estimates (23) and (24),
as well as the very strange arguments to get the
conclusion on page SMC\_ 2016 001606.} Alecsa~-- \cite[Theorems 3.3--3.12]%
{Ale17}; Karakaya et al.~-- \cite[Theorem 2.4]{KaAtDoBo17}; Okeke~\&
Abbas~-- \cite[Proposition 2.1]{OkeAbb17}; Sharma~\& Imdad~--
\cite[Proposition 4.9]{ShaImd17}; Abass et al.~-- \cite[Remark
2]{AbMeMe18}; Alagoz et al.~-- \cite[Theorem 2.1]{AlGuAk18}; Do{\u
g}an~-- \cite[Theorem 3.3.1]{Dog18}; Fathollahi~\& Rezapour~--
\cite[Propositions 2.1--2.3, 3.1, Theorem 3.2]{FatRez18}; Garodia~\&
Uddin~-- \cite[Theorem 3.1]{GarUdd18}, \cite[Theorem 3.1]{GarUdd20}
and \cite[Theorem 3.1]{GarUdd20b}; Hussain et al.~-- \cite[Theorem
3.4]{HuUlAr18}; Kumar~\& Chauhan~-- \cite[Theorems 1, 2]{KumCha18};
Wahab~\& Rauf~-- \cite[Theorem 3.4]{WahRau18}; Yildirim~--
\cite[Theorem 2]{Yil18};\footnote{Notice the strange estimation $d(x_{n}%
,p)\leq b_{n}$ with $b_{n}$ mentioned in (2.16) using the first inequality in
(2.15); a similar remark is valid for the estimation $d(x_{n},p)\leq a_{n}$
with $a_{n}$ from (2.17). A similar observation is valid for the estimations
$d(x_{n},p)\leq c_{n}$, $d(x_{n},p)\leq b_{n}$ and $d(x_{n},p)\leq a_{n}$
mentioned in the proof of \cite[Theorem 2]{YilAbb18}.} Yildirim~\& Abbas~--
\cite[Theorem 2]{YilAbb18}; Akhtar~\& Khan~-- \cite[Theorem 3.1--3.3]%
{AkhKha19};\footnote{Observe that $h\in\lbrack0.1)$ in the definition of a
generalized $C^{q}$-mapping and in the proof of Theorem 3.1 one takes
\textquotedblleft$\lambda=\max\{h,\frac{h}{1-h}\}$", and so $\lambda
=h/(1-h)\geq1$ for $h\in\lbrack1/2,1)$. Consequently, the estimate (3.7) does
not ensure that $\lim_{n\rightarrow\infty}d(x_{n+1},p)=0$ and that
$\lim_{n\rightarrow\infty}\frac{a_{n}}{b_{n}}=0$. Even if not mentioned,
probably $\lambda$ is the same in the proofs of Theorems 3.2 and 3.3.}
Asaduzzaman~\& Ali~-- \cite[Theorem 3.3]{AsaAli19}; Atalan~-- \cite[Theorem
3.3]{Ata19}; Atalan~\& Karakaya~-- \cite[Theorem 2.3]{AtaKar19b}; Bhutia~\&
Tiwary~-- \cite[Theorem 2.2--2.5]{BhuTiw19}; Ert\"{u}rk~\& G\"{u}rsoy~--
\cite[Theorem 2.3]{ErtGur19}; G\"{u}rsoy et~al.~-- \cite[Theorem
6]{GuEkKhKa19}; Gutti~\& Gedala~-- \cite[Theorem 5.1]{GutGed19}; Kumar~\&
Chugh~-- \cite[Theorem 2.2]{KumChu19}; Malik~\& Choudhary~-- \cite[Theorem
6]{MalCho19}; Mebawondu~\& Mewomo~-- \cite[Theorem 3.2]{MebMew19} and
\cite[Theorem 3.5]{MebMew19b};\footnote{In \cite[p.\ 302]{MebMew19} one must
replace $n\rightarrow0$ by $n\rightarrow\infty$ two times, while in
\cite[p.\ 11]{MebMew19b} one must interchange $\rightarrow\infty$ and
$\rightarrow0$ three times.} Okeke~-- \cite[Theorem 3.3]{Oke19}; Piri~et
al.~\cite[Lemmas 3.1, 3.2, Theorem 3.3]{PiDaRaGh19}; Aibinu~\& Kim~--
\cite[Theorem 3.2]{AibKim20}; Alagoz et~al.~-- \cite[Theorem 3.1]{AlGuAk20};
Ali~\& Ali~-- \cite[Theorem 2.3]{AliAli20} and \cite[Theorem 2.4]{AliAli20b};
Atalan~\& Karakaya~-- \cite[Theorem 2.4]{AtaKar20}; Chairatsiripong et~al.~--
\cite[Theorem 3.1]{ChYaTh20}; Deshmukh~-- \cite[Theorem 4.16]{Des20}; Dogan~--
\cite[Theorem 4]{Dog20}; Garodia et~al.~-- \cite[Theorem 3.1]{GaUdKh20};
G\"{u}rsoy et~al.~-- \cite[Theorem 2.3]{GuKhErKa20}; Ofem~\& Igbokwe~--
\cite[Theorem 3.1]{OfeIgb20} and \cite[Theorem 3.2]{OfeIgb21}; Sharma
et~al.~-- \cite[Theorems 2, 3]{ShMiMiPa20}; Shatanawi et~al.~-- \cite[Theorem
1]{ShBaTa20}; Udofia~\& Igbokwe~-- \cite[Theorem 3.1]{UdoIgb20}, \cite[Theorem
5.1]{UdoIgb21}, \cite[Theorem 3.1]{UdoIgb22}, \cite[Theorem 4.1]{UdoIgb23}~\&
\cite[Theorem 3.1]{UdoIgb23b}; Ali et~al.~-- \cite[Theorem 2.3]{AlAlRo21}; Ali
et~al.~-- \cite[Theorem 3.1]{AlAlUd21}; Bantaojai et~al.~-- \cite[Theorem
3.1]{BaGaUdPaYi21}; Hac{\i}o{\u{g}}lu~-- \cite[Theorem 7]{Hac21}; Hussain
et~al.~-- \cite[Theorem 9]{HuHuAl21}; Jubair et~al.~-- \cite[Theorem
3.2]{JuAlAl21}; Kalsoom et~al.~-- \cite[Theorem 19]{KaRaSuBiGhInAl21};
Maibed~\& Thajil~-- \cite[Theorems 2.8, 2.9]{MaiTha21}; Thajil~\& Maibed~--
\cite[Theorems 2.1, 2.2]{ThaMai21} and \cite[Theorems 2.3--2.5]{ThaMai21b};
Udofia et~al.~-- \cite[Theorem 10]{UdOfIg21b}; Abbas et~al.~-- \cite[Theorem
2]{AbAsDe22}; Ali et~al.~-- \cite[Theorems 2.2--2.5]{AlHuKaCh22}; Ali
et~al.~-- \cite[Theorem 3.3]{AlJuAl22};\footnote{Notice the interesting
expression: $\left\Vert \tau_{2,n}-t\right\Vert \leq\delta^{2(n+1)}\left(
1-(1-\delta\right)  \theta_{n}\mu_{n})^{n+1}\left\Vert \tau_{2,0}-t\right\Vert
=\alpha_{2,n}$, $n\in\mathbb{Z}_{+}$!} Atalan~\& Kilic~-- \cite[Theorem
3]{AtaKil22}; Botmart et~al.~-- \cite[Theorem 4.3]{BoShBaRe23}; Celik~\&
Simsek~-- \cite[Theorems 2.3]{CelSim22}; Hammad et~al.~-- \cite[Theorem
2]{HaReDe22}; Jia et~al.~-- \cite[Proposition 10]{JiShAhShBo22}; Jubair
et~al.~-- \cite[Theorem 16]{JuAlKu22}; Kaur~\& Chandok~-- \cite[Theorem
4]{KauCha22b}; Kaur~\& Chauhan~-- \cite[Case 1]{KauCha22}; Khan~--
\cite[Theorem 2]{Kha22}; Khan et~al.~-- \cite[Theorem 3]{KhAkDiSh22}; Kim
et~al.~-- \cite[Theorem 2.2]{KiDaPaSa22}; Maibed~\& Hussein~-- \cite[Theorem
(2.17)]{MaiHus22}; Ofem et~al.~-- \cite[Theorem 3.2]{OfIsAlAh22}; Rahimi
et~al.~-- \cite[Theorems 3.1, 3.3]{RaReDaGh22} and \cite[Theorem
2]{RaReDaGh24}; Sahu et~al.~-- \cite[Theorem 3.3]{SaBaGu22}; Salem~\&
Maibed~-- \cite[Theorem 2.5]{SalMai22}; Srivastava~-- \cite[Proposition
3.1]{Sri22}; Akram~-- \cite[Theorem 2.1]{Akr23};\footnote{It is worth
mentioning the following definition, where [40] is our reference \cite{Ber04}:
\textquotedblleft Definition 2.2 ([40]). Let $\left\{  p_{n}\right\}
_{n=0}^{\infty}$ and $\left\{  a_{n}\right\}  _{n=0}^{\infty}$ be two
\emph{real sequences} with $\lim_{n\rightarrow\infty}p_{n}=t^{\ast}$ and
$\lim_{n\rightarrow\infty}a_{n}=t^{\ast}$. If $\left\{  \mu_{n}\right\}  $ and
$\left\{  \nu_{n}\right\}  $ are two positive sequences converging to $0$
satisfying $\left\Vert p_{n}-t^{\ast}\right\Vert \leq\mu_{n}$ and $\left\Vert
a_{n}-t^{\ast}\right\Vert \leq\nu_{n}$, $\forall n\in\mathbb{N}$. Then
$\left\{  p_{n}\right\}  $ converges to $t^{\ast}$ faster than $\left\{
a_{n}\right\}  $ if $\left\{  \mu_{n}\right\}  $ converges faster than
$\left\{  \nu_{n}\right\}  $." See also \cite[Def.\ 1.1]{BhuTiw19}.} Akewe
et~al.~-- \cite[Sec.\ 4.4, Cases 1--9]{AkSaFa23}; Ali et~al.~-- \cite[Theorems
2--5]{AlAlDaTuZaMa23}; Anku et~al.~-- \cite[Theorem 7]{AnNaKa23}; Ansar~\&
Mas'ud~-- \cite[Theorem 3]{AnsMas23}; Atalan~\& Erba\c{s}~-- \cite[Theorem
3.5]{AtaErb23}; Botmart et~al.~-- \cite[Theorem 4.3]{BoShBaRe23};
Chairatsiripong et~al.~-- \cite[Theorem 4.6]{ChYaTh23}; Chauhan et~al.~--
\cite[Sec.\ 2.1]{ChKuImAs23}; Deshmukh et~al.~-- \cite[Theorems 3.21,
4.15]{DeGoRa23}; Fan~\& Wang~-- \cite[Theorem 3.1, Remark 3.1]{FanWan23};
Gautam~\& Kaur~-- \cite[Theorem 3]{GauKau23}; Gundogdu~-- \cite[Theorem
3.2.23]{Gun23}; Hammad~\& Kattan~-- \cite[Theorem 2]{HamKat23b} and
\cite[Theorem 2]{HamKat23}; Joodi~\& Maibed~-- \cite[Theorem 2.6]{JooMai23}
and \cite[Theorem 3.15]{JooMai23b}; Khan et~al.~-- \cite[Theorem
3.1]{KhAhAb23}; Maibed~\& AL-Hameedwi~-- \cite[Theorem 2.10]{MaiHam23}; Maibed
et~al.~-- \cite[Theorem 2.2]{MaJoSm23}; Maibed~\& Salem~-- \cite[Theorems
2.14--2.16]{MaiSal23}; Mebawondu et~al.~-- \cite[Theorem 3.3]{MePiNaOnAd23};
Ofem~-- \cite[Theorem 3.1]{Ofe23}; Okeke et~al.~-- \cite[Theorem
3.5]{OkUdAlAl23} and \cite[Theorem 3]{OkUdAlAl24}; Panigrahy~\& Mishra~--
\cite[Theorems 3.4, 3.6]{PanMis23}; Qawasmeh et~al.~-- \cite[Theorem
4]{QaBaBaQaHa23}; Sahu~\& Banerjee~-- \cite[Lemmas 3.1, 3.3, Theorem
3.5]{SahBan23} and \cite[Theorem 4]{SahBan24}; Salman~\& Abed~-- \cite[Theorem
4.1]{SalAbe23}; Sharma et~al.~-- \cite[Theorem 4]{ShRaBeKa23}; Ahmad et~al.~--
\cite[Theorem 3]{AhUlAhAlMl24}; Ahmed~-- \cite[Proposition 4.1]{Ahm24};
Albaqeri et~al.~-- \cite[Theorem 3.2]{AlHaReSe24}; Atalan et~al.~--
\cite[Theorem 11]{AtHaErGuMi24}; Begum~-- \cite[Theorems 3.1.1, 4.0.1]{Beg24};
El~Harmouchi et~al.~-- \cite[Theorem 3.3]{HaOuChMo24}; Gautam~\& Vineet~--
\cite[Theorem 2.3]{GauVin24};\footnote{In this paper one introduced the
following \textquotedblleft Definition 1.6 Let $\{p_{m}\}$ and $\{q_{m}\}$ be
two sequences in a Banach Space $X$ such that both $\{p_{m}\}$ and $\{q_{m}\}$
converge to the same point $p$. We say that $\{p_{m}\}$ converges to $p$
faster than $\{q_{m}\}$ if, for any positive real number $\epsilon_{2}%
>0$\medskip, there exists $\epsilon_{1}>0$ and $a\in\mathbb{N}$ such that
$\epsilon_{1}<\epsilon_{2}$, $\left\Vert p_{m}-p\right\Vert <\epsilon_{1}$,
and $\left\Vert q_{m}-p\right\Vert <\epsilon_{2}$ for all $m\geq a$".
Moreover, one \textquotedblleft demonstrate that the definition 1.6 is
consistent with the definition 1.5", where definitions 1.4 and 1.5 are
equivalent with \cite[Definitions 2.5 and 2.7]{Ber04} (in which
\textquotedblleft are available" is replaced by \textquotedblleft(best ones
available)", but without \textquotedblleft(converging to zero)"). Indeed,
having \textquotedblleft$\{p_{m}\}$ and $\{q_{m}\}$ ... two sequences in ...
$X$ such that both $\{p_{m}\}$ and $\{q_{m}\}$ converge to the same point
$p$", it is easy to prove that \textquotedblleft$\{p_{m}\}$ converges to $p$
faster than $\{q_{m}\}$" using \cite[Definition 1.6]{GauVin24}: Take
$\epsilon_{2}>0$; because $q_{m}\rightarrow p,$ there exists $a_{2}%
\in\mathbb{N}$ such that $\left\Vert q_{m}-p\right\Vert <\epsilon_{2}$ for
$m\geq a_{2}$. Consider $\epsilon_{1}:=\epsilon_{2}/2$ $(<\epsilon_{2})$;
because $p_{m}\rightarrow p$, there exists $a_{1}\in\mathbb{N}$ such that
$\left\Vert p_{m}-p\right\Vert <\epsilon_{1}$ for $m\geq a_{1}$; setting
$a:=\max\{a_{1},a_{2}\},$ the conclusion follows. Also notice that in the
proof of Theorem 2.3 one uses definition 1.5.} Keten~Copur et~al.~--
\cite[Theorem 3.3]{KeHaGu24}; Murali~\& Muthunagai~-- \cite[Theorems 3.3,
3.5]{MurMut24}; Okeke et~al.~-- \cite[Theorem 3.2]{OkAnUdOl24}; Panwar~\&
Bhokal~-- \cite[Theorem 2.4]{PanBho24}; Rani et~al.~-- \cite[Theorem
3.3]{RaKaBh24}; Agwu et~al.~-- \cite[Proposition 3.1]{AgSaIs25}; Akram~\&
Ahmad~-- \cite[Theorem 2.1]{AkrAhm25}; Alam~\& Rohen~-- \cite[proof of Theorem
6]{AlaRoh25}; Alam et~al.~-- \cite[Theorem 3.5]{AlRoTo25}; Chairatsiripong
et~al.~-- \cite[Theorem 2.3]{ChYaPaTh25}; Filali et~al.~-- \cite[Theorem
4]{FiElAlAlAlKh25};\footnote{In \cite{FiElAlAlAlKh25} one finds the next
definition in which [2] is our reference \cite{Ber04}: \textquotedblleft
Definition 1 ([2]). Consider two iteration sequences, $\left\{  u_{n}\right\}
$ and $\left\{  v_{n}\right\}  $, that both converge to the same point,
$u^{\ast}$. If there exist two real-valued sequences, $\left\{  \zeta
_{n}\right\}  $ and $\left\{  \eta_{n}\right\}  $, such that $\left\Vert
u_{n}-u^{\ast}\right\Vert \leq\zeta_{n}$ and $\left\Vert v_{n}-u^{\ast
}\right\Vert \leq\eta_{n}$ for all $n=1,2,3,...$, then the sequence $\left\{
u_{n}\right\}  $ is said to converge more rapidly than $\left\{
v_{n}\right\}  $ if $\lim_{n\rightarrow\infty}\frac{\zeta_{n}}{\eta_{n}}=0$."
As seen in the proof of our Proposition \ref{p-fast}, having \textquotedblleft
two iteration sequences, $\left\{  u_{n}\right\}  $ and $\left\{
v_{n}\right\}  $, that both converge to the same point, $u^{\ast}$",
\textquotedblleft there exist two real-valued sequences, $\left\{  \zeta
_{n}\right\}  $ and $\left\{  \eta_{n}\right\}  $" with the mentioned
properties (even convergent to $0$), and so $\left\{  u_{n}\right\}  $
converge more rapidly than $\left\{  v_{n}\right\}  $.} Ishtiaq et~al.~--
\cite[Theorem 4.5]{IsBaHuAl25}; Mary~\& Uthayakumar~-- \cite[Theorem
3.3]{MarUth25}; Nawaz et~al.~-- \cite[Theorem 4.1]{NaGdUl25} and \cite[Theorem
5]{NaUlGd25b}; Okeke et~al.~-- \cite[Theorem 4]{OkAlAn25}; Sharma et~al.~--
\cite[Theorem 3]{ShRaKaBeRa25};

\smallskip As mentioned in Section \ref{sec1}, Dogan~\& Karakaya obtained that
\textquotedblleft the iteration schemes $\{k_{n}\}_{n=0}^{\infty}$ and
$\{l_{n}\}_{n=0}^{\infty}$ have the same rate of convergence to $p$ of $\wp$"
in \cite[Theorem 2.4]{DogKar18} because they found the same upper estimates
for $\left\Vert k_{n+1}-p\right\Vert $ and $\left\Vert l_{n+1}-p\right\Vert $
when $l_{0}=k_{0}$ (see \cite[page 156]{DogKar18}); the same argument is used
for getting the same conclusion by Kumar~\& Chauhan (Gonder) (see \cite[page
947]{KumCha20}).

\medskip

It is worth repeating that Popescu (in \cite{Pop07}) recalls \cite[Definition
2.7]{Ber04}, mentions its inconsistency, introduces his direct comparison of
iterative processes in \cite[Definition 3.5]{Pop07}, and uses this definition
in \cite[Theorem 3.7]{Pop07}.

\medskip Other papers in which \cite[Definition 3.5]{Pop07} is used, possibly
without citing it (but possibly recalling \cite[Definition 2.5 or/and
Definition 2.7]{Ber04}), are: Xue~-- \cite[Theorems 2.1, 2.2]{Xue08};
Rhodes~\& Xue~-- \cite[Theorems 2.1, 2.2, 3.1, 3.2]{RhoXue10}; Thong~--
\cite[Theorems 2.1, 2.3, 2.5]{Tho12}; Alotaibi et al.~-- \cite[Theorem
3.1]{AlKuHu13}; Hussain et al.~-- \cite[Theorems 14--17]{HuKuKu13}%
\footnote{Note the strange quantity $\big\Vert\frac{\text{JN}_{n+1}%
-p}{\text{JI}_{n+1}-p}\big\Vert$, the numerator and denominator being in
$(X,\left\Vert \cdot\right\Vert )$ \textquotedblleft an arbitrary Banach
space\textquotedblright.}; Phuengrattana~\& Suantai~-- \cite[Theorems 2.4,
2.6]{PhuSua13}; Khan et al.~-- \cite[Theorem 3.1]{KhKuHu14}; Fukhar-ud-din~\&
Berinde~-- \cite[Theorems 2.5, 2.7]{FukBer16}; G\"{u}rsoy~-- \cite[Theorem
2.4]{Gur16}; Khan et al.~-- \cite[Theorem 3]{KhGuKa16}; G\"{u}rsoy et al.~--
\cite[Theorem 2.3]{GuKhFu17}; Kosol~-- \cite[Theorem 2.2]{Kos18};\footnote{In
\cite{Ber04} one finds \textquotedblleft$\left\Vert u_{n}-p\right\Vert
\leq\left\Vert v_{n}-p\right\Vert ,\ \forall n$. \ (2.4)", meaning that
$\{u_{n}\}_{n=0}^{\infty}$ converges better than $\{v_{n}\}_{n=0}^{\infty}$ in
the sense of Rhoades \cite{Rho76}; in \cite{Ber04} there are only upper
estimates for $\left\Vert x_{n+1}-p\right\Vert $ and $\left\Vert
y_{n+1}-p\right\Vert $.} Pansuwan~\& Sintunavarat~-- \cite[Theorem
3.7]{PanSin18}; Atalan~\& Karakaya~-- \cite[Theorem 3.3]{AtaKar19};
Chumpungam~\& Kettapun~-- \cite[Theorem 3.6]{ChuKet19}; Dung~\& Hieu~--
\cite[Propositions 3.5, 3.6]{DunHie20}; Ert\"{u}rk~\& G\"{u}rsoy~--
\cite[Theorem 2.3]{ErtGur19}; Kumam et al.~-- \cite[Theorems 3.4,
3.5]{KuKhKu19}; G\"{u}rsoy et al.~-- \cite[Theorem 4]{GuErAb20}; Khai et
al.~-- \cite[Theorems 3.2, 3.8--3.10]{KhPhHuNgHi25}; Udomene \cite[Theorems
2.1, 2.2]{Udo25}.

\medskip

It is also worth noticing that by taking simple examples in
$\mathbb{R}$, Rafiq et~al.~-- \cite[Example 11]{RaLeDaDj10}; Hussain
et al.~-- \cite[Example 9]{HuRaDaLa11}; Chugh et al.~--
\cite[Example 4.1]{ChKuKu12}; Hussain et al.~-- \cite[Examples 3.1,
3.2]{HuChKuRa12}; Kang et al.~-- \cite[Example 11]{KaCiRaAlKw13};
Karakaya et al.~-- \cite[Example 4]{KaDoGuEr13}; Kumar et al.~--
\cite[Example 9]{KuLaRaHu13}; Do\u{g}an~\& Karakaya~-- \cite[Example
10]{DogKar14}; Prasad~\& Goyal~-- \cite[Example 2.1]{PraGoy16};
Wahab~\& Rauf~-- \cite[Example 11, Remarks 12--17]{WahRau16};
Chauhan et~al.~-- \cite[Example 3.1]{ChUtImAh17}; Sintunavarat~--
\cite[Example 13]{Sin17};
Ullah~\& Arshad~-- \cite[Example 4.1]{UllArs17}, \cite[Example 4.1]%
{UllArs18}~\& \cite[Example 3.4]{UllArs18b}; Akeke~\& Eke~-- \cite[Example
3.1]{AkeEke18}; Akutsah et al.~-- \cite[Examples 1, 2]{AkNaAfMe21}; Tidke~\&
Patil~-- \cite[page 14]{TidPat23} and \cite[page 114]{TidPat23b}; Saif
et~al.~-- \cite[Example 3.1]{SaAlAlUd24}; Khan et~al.~-- \cite[Examples 2, 3,
5]{KhDiIqAb25}; Nawaz et~al.~-- \cite[Example 2]{NaUlGd25} and \cite[page
1956]{NaUlGd25b} \textquotedblleft prove" that certain iteration processes are
faster than other ones.

\medskip\textbf{Final remark}. We wish to point out that this paper is not
about the correctness of the results in the cited papers; we did not
check the proofs of the results. Our aim is to emphasize again, as
Popescu \cite{Pop07} and Phuengrattana--Suantai \cite{PhuSua13} did,
that Berinde's method is inconsistent, and so, what is obtained
using it, is useless from our point of view. The other remarks
mainly concern wrong attributions of notions as well as the fact
that one can not claim the validness of general assertions using
some examples; of course, (counter)examples are used to invalidate
results.

\end{document}